# FREE EXTREME VALUES

By Gerard Ben Arous and Dan Virgil Voiculescu[1]

*New York University and University of California at Berkeley*

Free probability analogs of the basics of extreme-value theory are obtained, based on Ando's spectral order. This includes classification of freely max-stable laws and their domains of attraction, using "free extremal convolutions" on the distributions. These laws coincide with the limit laws in the classical peaks-over-threshold approach. A free extremal projection-valued process over a measure-space is constructed, which is related to the free Poisson point process.

**1. Introduction.** Free probability theory [17, 18] is a highly noncommutative parallel to a large part of basic classical probability theory. The aim of this paper is to add a somewhat unexpected new entry to the classical-free dictionary: extreme-value theory.

Classical probability theory has constructed a complete picture for the behavior of extreme values of i.i.d. samples (see [7, 14]). This has been motivated by the huge importance of this question for statistics (catastrophic events, insurance, reliability and, more recently, finance and risk theory). The results provide an exhaustive description of the possible limits of normalized extreme values as being the extreme value distribution of their domain of attraction. In this note we develop a free analog to classical extreme-value theory. Also, surprisingly, the free extreme-value distributions turn out to have classical realizations already used by statisticians.

Addition and multiplication of free noncommutative random variables give rise to additive and multiplicative free convolution operations on probability measures on $\mathbb{R}$. With respect to the so-called *spectral order* [1] there are also analogs of the min and max operations for self-adjoint operators. Passing to the distributions, this gives rise in the case of free random variables to the "extremal free convolution" operation that we consider here.

Received February 2005; revised May 2005.
[1]Supported by NSF Grant 00-79945.

These turn out to be quite easy to express using distribution functions, without the analytic function machinery necessary for the other free convolutions. This can also be viewed as applying the limiting process that leads to idempotent analysis [8, 10] to additive free convolution.

The natural processes in this context are noncommutative extremal processes over an ordered set. In particular, we show there are some remarkable such processes indexed by the measurable sets of a measure space: the projection-valued and the triangular free extremal processes over a set. We show that the projection-valued free extremal process over a set can be realized via the free Poisson process and therefore has a natural asymptotic random matrix realization.

The analytical aspect of proving the analog of the classical results of Frechet, Fisher–Tippet and Gnedenko has as a main step establishing that the classical and free domains of attraction are actually the same. This is similar to the situation that occurs for the additive domains of attraction [3]. After finding the domains of attraction, obtaining limit laws with the same normalization constants as in the classical case is rather direct. The laws we find are called generalized Pareto laws in the statistics literature. More precisely, the Frechet law, the Gumbel law and the Weibull law have as free correspondents the Pareto law, the exponential law and the Beta law. Note that, reminiscent of the additive case [4], there is a certain shrinking of supports for the Gumbel and Weibull distributions.

These generalized Pareto distributions, which appear as the only free max-stable distributions, have played a role in statistics for a long time in the P.O.T. (peaks-over-threshold) approach to extreme value theory [2, 12].

After some preliminaries on spectral order in Section 2, we discuss extremal free convolution operations in Section 3. Then in Section 4 we introduce extremal noncommutative processes and in Section 5 we consider extremal free processes over a set. Section 6 is devoted to the study of the iteration of free extremal convolution and to the description of free max-stable distributions and their domain of attraction. In Section 7 we recall the classical P.O.T. context where the same laws appear.

**2. Preliminaries on $\wedge$ and $\vee$.** We shall work in a tracial $W^*$-probability space $(M, \tau)$, that is, $M$ is a von Neumann algebra with an ultraweakly continuous faithful tracial state $\tau$. For the basic von Neumann algebra facts we will use, the reader may consult [9] or [16]; for the basic free probability, references are [17] and [18]. It will be convenient to assume $M$ acts on a Hilbert space $\mathcal{H}$. For definiteness we may take $\mathcal{H} = L^2(M, \tau)$, which is naturally a left $M$-module since the $L^2$-space is the completion of $M$ viewed as a pre-Hilbert space with respect to the scalar product $\langle a, b \rangle = \tau(b^*a)$.

By $\text{Proj}(M)$ we shall denote the set of self-adjoint projections $P = P^2 = P^* \in M$. If $P, Q \in \text{Proj}(M)$, the self-adjoint projections defined by $(P \vee$ $Q)\mathcal{H} = \overline{P\mathcal{H} + Q\mathcal{H}}$ and $(P \wedge Q)\mathcal{H} = P\mathcal{H} \cap Q\mathcal{H}$ are in $\mathrm{Proj}(M)$. On $\mathrm{Proj}(M)$ the order relation $P \leq Q$, that is, $Q - P$ is a positive operator, coincides with $P\mathcal{H} \subset Q\mathcal{H}$, and $P \vee Q$ and $P \wedge Q$ are the max and min with respect to this order.

Let $M_h = \{m \in M | m = m^*\}$ denote the self-adjoint elements. If $a \in M_h$ and $\omega \subset \mathbb{R}$ is a Borel set, then $E(a;\omega) \in \mathrm{Proj}(M)$ will denote the corresponding spectral projection. By definition, the spectral order relation [1] on $M_h$ is $a \prec b$ if $E(a;[t,\infty)) \leq E(b;[t,\infty))$ for all $t \in \mathbb{R}$. Since $a = \int_0^\infty E(a;[t,\infty))\,dt$ if $a \geq 0$, it is easy to see that $a \prec b \Rightarrow a \leq b$ in case $a \geq 0$, $b \geq 0$ and, from here, also for general $a,b \in M_h$. Also $a \prec b \Rightarrow f(a) \prec f(b)$ if $f:\mathbb{R} \to \mathbb{R}$ is an increasing Borel function. All this also extends to self-adjoint unbounded operators affiliated with $M$.

The operations $\wedge$ and $\vee$ have a natural extension to $M_h$ (and even to affiliated self-adjoint operators). It is given in [1] in the case where $M$ is a matrix algebra, but the extension to the case of a general von Neumann algebra is quite straightforward and must have occurred to many people. If $a,b \in M_h$, then $a \wedge b$ is defined by $E(a \wedge b;[t,\infty)) = E(a;[t,\infty)) \wedge E(b;[t,\infty))$ for all $t \in \mathbb{R}$. To see that $a \wedge b$ is well defined, observe that the right-hand side is projection-valued, decreasing and left-continuous in the strong operator topology as a function of $t$. Similarly one defines $a \vee b$ by $E(a \vee b;(t,\infty)) = E(a;(t,\infty)) \vee E(b;(t,\infty))$. Here, checking that $a \vee b$ is well defined boils down to checking right continuity in $t$. Note that these definitions work in an arbitrary von Neumann algebra, that is, we did not use the tracial state $\tau$.

In the tracial context we also have $E(a \vee b;[t,\infty)) = E(a;[t,\infty)) \vee E(b;[t,\infty))$. Since $E(c;(t-\varepsilon,\infty)) \downarrow E(c;[t,\infty))$ as $\varepsilon \downarrow 0$, it is immediate that the left-hand side is greater than or equal the right-hand side, while the fact that this must be an equality follows from the following inequalities that involving the tracial state:

$$0 \leq \tau(E(a \vee b;(t-\varepsilon,\infty)) - E(a;[t,\infty)) \vee E(b;[t,\infty)))$$
$$\leq \tau(E(a;(t-\varepsilon,t)) + E(b;(t-\varepsilon,t))).$$

Using $(-a) \vee (-b) = -(a \wedge b)$ we also get $E(a \wedge b;(t,\infty)) = E(a;(t,\infty)) \wedge E(b;(t,\infty))$. Passing to orthocomplements also gives in the tracial case that

$$E(a \wedge b;(-\infty,t)) = E(a;(-\infty,t)) \vee E(b;(-\infty,t)),$$
$$E(a \wedge b;(-\infty,t]) = E(a;(-\infty,t]) \vee E(b;(-\infty,t]),$$
$$E(a \vee b;(-\infty,t]) = E(a;(-\infty,t]) \wedge E(b;(-\infty,t]),$$
$$A(a \vee b;(-\infty,t)) = E(a;(-\infty,t)) \wedge E(b;(-\infty,t)).$$

Note, that since only spectral projections are involved, the definitions and properties of $\wedge$ and $\vee$ also extend to unbounded self-adjoint operators affiliated with $(M,\tau)$. Note also that if $f:\mathbb{R} \to \mathbb{R}$ is an increasing Borel function, then $f(a) \wedge f(b) = f(a \wedge b)$ and $f(a) \wedge f(b) = f(a \vee b)$.

LEMMA 2.1. *If $P, Q \in \operatorname{Proj}(M)$ are freely independent in $(M, \tau)$, then*

$$\tau(P \vee Q) = \min(\tau(P) + \tau(Q), 1) \quad and \quad \tau(P \wedge Q) = \max(0, \tau(P) + \tau(Q) - 1).$$

This is not a new result. It can be obtained using free convolution, additive or multiplicative. Indeed, $P \wedge Q = E(P + Q, \{2\}) = E(PQP, \{1\})$. Since the Cauchy transforms of the distributions of $PQP$ and $P + Q$ are algebraic, $\tau(P \wedge Q)$ is given by the residue at 1 or, respectively, at 2 of the corresponding Cauchy transform (see 3.4.1 in [17] or 3.6.7 in [18]).

The free independence of two projections of given trace completely determines the trace on the algebra they generate. Since in this paper we focus on $P \wedge Q$ and $P \vee Q$, it will be convenient to also have at hand a more relaxed concept.

DEFINITION 2.2. Two projections $P, Q \in \operatorname{Proj}(M)$ are in general position if the equivalent conditions $\tau(P \vee Q) = \min(\tau(P) + \tau(Q), 1)$ and $\tau(P \wedge Q) = \max(0, \tau(P) + \tau(Q) - 1)$ are satisfied. We will also say that two unbounded self-adjoint operators $a, b$ affiliated with $M$ are in general spectral position if, for each $t \in \mathbb{R}$, the projections $E(a; [t, \infty))$ and $E(b; [t, \infty))$ are in general position.

The preceding definition for self-adjoint operators may seem to depend on choosing $[t, \infty)$ instead of $(t, \infty)$, but this is actually inessential.

LEMMA 2.3. *If the unbounded self-adjoint operators $a, b$ affiliated with $(M, \tau)$ are in general position, then $E(a; (t, \infty))$ and $E(b; (t, \infty))$ are in general position. In particular, $-a$ and $-b$ are in general position.*

The preceding lemma follows from the inequalities

$$0 \leq \tau(E(a; (t, \infty)) \vee E(b; (t, \infty)) - E(a; [t + \varepsilon, \infty)) \vee E(b; [t + \varepsilon, \infty)))$$
$$\leq \tau(E(a; (t, t + \varepsilon)) + E(b; (t, t + \varepsilon))).$$

Returning to the spectral order we have the following lemma.

LEMMA 2.4. *If $a, b \in M$, $a \geq 0$, $b \geq 0$, then $(2^{-1}(a^p + b^p))^{1/p} \uparrow a \vee b$ as $p \to +\infty$.*

With the superfluous condition that $M$ be finite-dimensional, this is Lemma 6.15 in [1]. The proof in [1] that $2^{-1/p}(a^p + b^p)^{1/p}$ is increasing, works in general and the rest of the proof also works with minor adjustments. If $X$ denotes the limit, then $2^{-1}(a^p + b^p) \leq (a \vee b)^p$ gives $X \leq a \vee b$ since $t \to t^{1/p}$ is operator-increasing for $p \geq 1$. On the other hand, if $k \geq 1$, $X^k = s - \lim_{p \to \infty} 2^{-k/p}(a^p + b^p)^{k/p} = s - \lim_{p \to \infty} 2^{-1/p}(a^{kp} + b^{kp})^{1/p} \geq a^k$, so that if $\xi \in E(X;(-\infty,t])\mathcal{H}$, then $t^k\|\xi\|^2 \geq \langle X^k\xi,\xi\rangle \geq \langle a^k\xi,\xi\rangle$ for all $k \in \mathbb{N}$, so that $\xi \in E(a;(-\infty,t])\mathcal{H}$. Similarly, $\xi \in E(b;(-\infty,t])\mathcal{H}$, so that $E(X;(-\infty,t]) \leq E(a;(-\infty,t]) \wedge E(b;(-\infty,t])$. This easily gives $X \geq a \vee b$.

COROLLARY 2.5. *If $a,b \in M_h$, then*
$$a \vee b = s - \lim_{p\to\infty} p^{-1}\log(\exp(pa) + \exp(pb)).$$

Indeed this follows from $\log(\exp(a \vee b)) = a \vee b$, $2^{-1/p} \to 1$ and the preceding lemma applied to $\exp a$, $\exp b$.

**3. Extremal free convolutions.** By $\mathrm{Prob}(\mathbb{R})$ we shall denote the probability measures on $\mathbb{R}$ and by $\mathrm{Prob}_c(\mathbb{R})$ denote those with compact support. If $\mu \in \mathrm{Prob}(\mathbb{R})$, then $F(t) = \mu((-\infty,t])$ is its distribution function.

DEFINITION 3.1. *If $\mu,\nu \in \mathrm{Prob}(\mathbb{R})$ have distribution functions $F, G$, then the upper and the lower extremal free convolutions $\mu \boxvee \nu$ and, respectively, $\mu \boxwedge \nu$ are given by the distribution functions $H(t) = \max(0, F(t) + G(t) - 1)$ and, respectively, $K(t) = \min(F(t) + G(t), 1)$.*

PROPOSITION 3.2. *If $a,b$ are self-adjoint unbounded operators affiliated with $(M,\tau)$ that are in general spectral position and if $\mu_a, \mu_b \in \mathrm{Prob}(\mathbb{R})$ are their distributions, then $\mu_a \boxvee \mu_b$ and $\mu_a \boxwedge \mu_b$ are the distributions of $a \vee b$ and, respectively, $a \wedge b$.*

The preceding proposition is immediate from the definitions. Using Lemma 2.1, Definition 2.2 and the adaptation of free independence to unbounded operators affiliated with $(M,\tau)$ [5] we have the following consequence.

COROLLARY 3.3. *If $a,b$ are freely independent unbounded self-adjoint operators affiliated with $(M,\tau)$, then the distributions of $a \vee b$ and $a \wedge b$ are equal to $\mu_a \boxvee \mu_b$ and, respectively, to $\mu_a \boxwedge \mu_b$ where $\mu_a, \mu_b$ are the distributions of $a, b$.*

The extremal free convolution operations also have natural descriptions without invoking the distribution functions.

PROPOSITION 3.4. *Let $\mu,\nu \in \mathrm{Prob}(\mathbb{R})$ and let $t = \inf\{x \in \mathbb{R} | (\mu + \nu) \times ((x,\infty)) \leq 1\}$ and $s = \sup\{y \in \mathbb{R} | (\mu + \nu)((-\infty,y)) \leq 1\}$. Then*
$$\mu \boxvee \nu = \chi_{(t,\infty)} \cdot (\mu + \nu) + (1 - (\mu + \nu)((t,\infty)))\delta_t$$

*and*
$$\mu \boxwedge \nu = \chi_{(-\infty,s)} \cdot (\mu + \nu) + (1 - (\mu + \nu)((-\infty,t)))\delta_s.$$

PROOF. Given that the map $x \to -x$ interchanges the two operations, it will suffice to check the first equality. If $F$, $G$ and $H$ are the distribution functions of $\mu$, $\nu$ and the right-hand side of the equality we are checking, respectively, it is immediate that

$$1 - H(x) = \min((1 - F(x)) + (1 - G(x)), 1),$$

which gives

$$H(x) = \max(F(x) + G(x) - 1, 0). \qquad \square$$

Translating Corollary 2.5 into free convolutions language, we obtain the following result.

PROPOSITION 3.5. *Let $\mu, \nu \in \text{Prob}_c(\mathbb{R})$. Then $\mu \boxdot \nu$ is the weak limit of*

$$(k^{-1} \log)_*((\exp k\bullet)_*\mu \boxplus (\exp k\bullet)_*\nu)$$

*as $k \to +\infty$.*

Note that the preceding proposition can be given additional precision using Lemma 2.4 with the $2^{-1/p}$ factor.

REMARK 3.6. It is interesting to note that in the classical context [7, 14], the operation on probability measures, analogous to $\boxdot$, corresponds to the multiplication of the distribution functions. Thus, viewing distribution functions as functions with values in $[0, 1]$, the passage from classical to free probability means replacing the multiplicative semigroup $([0, 1], \bullet)$ by the semigroup on $[0, 1]$ that results from the binary operation $(s, t) \to \max(0, s + t - 1)$.

## 4. Extremal noncommutative processes indexed by an ordered set.

DEFINITION 4.1. Let $(\mathcal{J}, \leq)$ be an ordered set and let $(M, \tau)$ be a tracial $W^*$-probability space. An $M$-valued upper extremal process indexed by $\mathcal{J}$ is a family $(Y(i))_{i \in \mathcal{J}}$ of self-adjoint unbounded operators affiliated with $M$ such that if $\alpha \in \mathcal{J}$ is the least upper bound of a set $\omega \subset \mathcal{J}$, then

$$Y(\alpha) = \vee \{Y(i) | i \in \omega\}$$

[for sets, $\vee$ is defined like for pairs on the spectral projections $E(\cdot; (t, \infty))$]. The process is projection-valued if $Y(i) \in \text{Proj}(M)$, $i \in \mathcal{J}$.

REMARK 4.2. The extension of the supremum $\vee$ to infinite families no longer has the property that spectral projections for $(t, \infty)$ may be replaced by spectral projections for $[t, \infty)$. For instance, $\bigvee_{n \in \mathbb{N}} (1 - 1/n)I = I$, but $E((1 - 1/n)I; [1, \infty)) = 0$ while $E(I; [1, \infty)) = I$.

REMARK 4.3. The preceding definition implies that if $\alpha, \beta \in \mathcal{J}$, $\alpha \leq \beta$, then $Y(\alpha) \prec Y(\beta)$ with respect to the spectral order. Indeed, it suffices to apply the definition to $\omega = \{\alpha, \beta\}$.

If $Y_1, Y_2$ are $M$-valued upper extremal processes indexed by the same ordered set $\mathcal{J}$, then we may define $Y_1 \vee Y_2$ by $(Y_1 \vee Y_2)(\alpha) = Y_1(\alpha) \vee Y_2(\alpha)$, which is also an upper extremal process indexed by $\mathcal{J}$. Clearly the construction extends to families of processes.

If $Y_j$, $j = 1, 2$, are upper extremal processes indexed by $\mathcal{J}$ with values in von Neumann algebras $(M_j, \tau_j)$, $j = 1, 2$, we may construct their *free sup* as the process $\sigma_1(Y_1) \vee \sigma_2(Y_2)$, where $\sigma_j : M_j \to M_1 * M_2$ are the canonical inclusions into $(M_1 * M_2, \tau_1 * \tau_2)$. Again, more generally we may take the free sup of a family of processes.

Consistent with these considerations is the spectral order on processes $Y_1 \prec Y_2$ defined as $Y_1(\alpha) \prec Y_2(\alpha)$ for all $\alpha \in \mathcal{J}$.

REMARK 4.4. If $Y$ is an $M$-valued upper extremal process indexed by $\mathcal{J}$, then $\mathcal{J} \ni \alpha \to E(Y(\alpha); (t, \infty)) \in \mathrm{Proj}(M)$ is a projection-valued upper extremal process indexed by $\mathcal{J}$.

EXAMPLE 4.5 (Spectral measures). If $a$ is an unbounded self-adjoint operator affiliated with $(M, \tau)$, then

$$\mathbb{R} \ni t \to E(a; (-\infty, t)) \in \mathrm{Proj}(M)$$

is a projection-valued upper extremal process indexed by $\mathbb{R}$.

REMARK 4.6. Note that $a \prec b$ is equivalent to $E(a; (-\infty, t)) \leq E(b; (-\infty, t))$, $t \in \mathbb{R}$. Thus, Ando's definition of the spectral order [1] can be interpreted as replacing the operator by the projection-valued process of Example 4.5 and using the natural order on such processes derived from the order on projections.

REMARK 4.7. The definition of upper extremal processes can be extended to processes where the $Y(\alpha)$ are self-adjoint unbounded operators affiliated with $M$. Note also that if $f : \mathbb{R} \to \mathbb{R}$ is an increasing function that is lower semicontinuous, then

$$\mathcal{J} \ni \alpha \to f(Y(\alpha))$$

transforms an upper extremal process into another upper extremal process. We denote this process by $f_*Y$.

**5. Free extremal processes over a set.** Let $(\mathcal{X}, \mathcal{B}, \mu)$ be a measure space with $\mu$ a positive $\sigma$-finite measure. The *free projection-valued upper extremal process* over $(\mathcal{X}, \mathcal{B}, \mu)$ will be a projection-valued upper extremal process $\mathcal{B}/\sim \ni \omega \to Y(\omega)$, where $\mathcal{B}/\sim$ denotes $\mathcal{B}$ modulo null sets with values in some $(M, \tau)$ so that if $\{\omega_k\}_{k \in K}$ are pairwise disjoint (modulo null sets), then the $Y(\omega_k)$, $k \in K$, are freely independent and $\tau(Y(\omega)) = \min(\mu(\omega), 1)$.

It is easy to see that it is sufficient to construct such a process when $\mathcal{X}$ is a finite set; the general case then is roughly the result of viewing $\mathcal{X}$ as arising from an inverse limit using finite partitions, while $M$ would be constructed as a direct limit (we leave the details to the reader). If $\mathcal{X}$ is finite and $\mathcal{X} = \{x_1, \ldots, x_n\}$, $\{x_j\} \in \mathcal{B}$, then let $P_1, \ldots, P_n \in \mathrm{Proj}(M)$ be freely independent with $\tau(P_j) = \min(\mu(\{x_j\}), 1)$. We then define, for $\omega \subset \mathcal{X}$,

$$Y(\omega) = \bigvee_{\{j \mid x_j \in \omega\}} P_j.$$

That $\tau(Y(\omega)) = \min(\mu(\omega), 1)$ is then a consequence of Lemma 2.1.

It is easily seen that under the additional requirement that $M$ be generated by $\{Y(\omega) | \omega \in \mathcal{B}\}$, the free projection-valued upper extremal process over $(\mathcal{X}, \mathcal{B}, \mu)$ is unique up to isomorphism.

The free projection-valued upper extremal process over $(\mathcal{X}, \mathcal{B}, \mu)$ is related to the free Poisson process over $(\mathcal{X}, \mathcal{B}, \mu)$. To explain this, we start by recalling some facts about the free Poisson process. Based on realizations of free Poisson variables using semicircular or circular elements (see [11], Remark 1.7, Corollary 1.8 in the main text and Lemma 1.4, Remark 1.5, Theorem 1.6 in the Appendix), there are constructions of free Poisson processes over a set (see Section 6.2 and 1° in Section 6.6 of [17]). We will use the construction that involves circular elements, which we summarize in the following theorem.

THEOREM 5.1 ([11, 17]). *Let $(C_\iota)_{\iota \in I}$ be $*$-freely independent circular variables and let $(\Omega_\iota)_{\iota \in I}$ be spaces of events with $\sigma$-algebras $\Sigma_\iota$ and probability measures $\nu_\iota$. Let further $(\Omega, \Sigma, \mu)$ be the disjoint union of $(\Omega_\iota, \Sigma_\iota, \nu_\iota)$. Let $(A_\iota, \tau_\iota)$ denote the $W^*$-probability space $L^\infty(\Omega_\iota, \Sigma_\iota, \nu_\iota)$ with $\tau_\iota$ the expectation functional given by $\nu_\iota$. Assume $(A_\iota)_{\iota \in I}$ and $(\{C_\iota, C_\iota^*\})_{\iota \in I}$ are freely independent and contained in $(M, \tau)$. If $\alpha \in \Sigma$, $\nu(\alpha) < \infty$, let*

$$\Pi(\alpha) = \sum_{\iota \in I} C_\iota \chi_{\alpha \cap \Omega_\iota} C_\iota^*,$$

*where $\chi_{\alpha \cap \Omega_\iota} \in A_\iota$ is the indicator function of $\alpha \cap \Omega_\iota$. Then:*

(i) $\Pi(\alpha)$ *has free Poisson distribution with parameters $a = \nu(\alpha)$, $b = 1$ in $(M, \tau)$.*

(ii) *If $\alpha$ is the disjoint union of $\alpha_k \in \Sigma$, $\nu(\alpha) < \infty$, then $\Pi(\alpha) = \sum_{k \in \mathbb{N}} \Pi(\alpha_k)$.*

(iii) *If $\alpha_k \in \Sigma$, $\nu(\alpha_k) < \infty$ are pairwise disjoint, $k \in \mathbb{N}$, then $(\Pi(\alpha_k))_{k \in \mathbb{N}}$ are freely independent.*

Using the facts on freeness of Gaussian and deterministic diagonal matrices, it was also noted (see [17], the end of Section 7.3) that one obtains asymptotic random matrix realizations of the free Poisson processes.

To avoid complicating notation, we gave in Theorem 5.1 a construction of the process for $\Omega$, which can be represented as a disjoint union of probability measure spaces. Clearly if $(\mathcal{X}, \mathcal{B}, \mu)$ is a $\sigma$-finite measure space, we can realize the free Poisson process over $(\mathcal{X}, \mathcal{B}, \mu)$ by finding $\Omega$ as above so that $\Omega \supset \mathcal{X}$, $\mathcal{B} \subset \Sigma$ and $\nu|\mathcal{X}$ extends $\mu$.

Let $(\mathcal{X}, \mathcal{B}, \mu)$ be a $\sigma$-finite measure space and let $\Pi : \mathcal{B}_f \to M_h$, where $\mathcal{B}_f$ are the sets in $\mathcal{B}$ with finite measure, be the free Poisson process, that is, $\Pi(\alpha)$ is a free Poisson noncommutative random variable with parameters $\mu(\alpha)$ and 1, and for disjoint $\alpha_1, \alpha_2, \ldots$ so that $\sum_k \mu(\alpha_k) < \infty$, we have $\sum_k \Pi(\alpha_k) = \Pi(\bigcup_k \alpha_k)$ and the $\Pi(\alpha_k)$ are freely independent in $(M, \tau)$. The following theorem gives the connection between the two kinds of processes.

THEOREM 5.2. *If $Y(\alpha)$ denotes the range projection of $\Pi(\alpha)$, then $Y(\alpha)$ is the free projection-valued upper extremal process over $(\mathcal{X}, \mathcal{B}, \mu)$. Moreover, if $\mu$ has no atoms, then the von Neumann algebras $\{\Pi(\alpha) | \alpha \in \mathcal{B}\}$ and $\{Y(\alpha) | \alpha \in \mathcal{B}\}$ are equal.*

PROOF. In view of the formula for the distribution of a free Poisson random variable (Section 2.7 in [17]), we have $\tau(Y(\alpha)) = \min(\mu(\alpha), 1)$. Also, clearly $\Pi(\alpha)$ depends on $\alpha$ only up to null sets and if $\alpha_1, \alpha_2, \ldots$ are disjoint, then $Y(\alpha_1), Y(\alpha_2), \ldots$ are freely independent, since $\Pi(\alpha_1), \Pi(\alpha_2), \ldots$ are freely independent. Further, if $\alpha_1, \alpha_2, \ldots$ are disjoint and $\alpha = \bigcup_k \alpha_k \in \mathcal{B}_f$, then $Y(\alpha) = \bigvee_k Y(\alpha_k)$. Indeed, $Y(\alpha)$ is the projection onto $(\ker \Pi(\alpha))^\perp$ and it suffices to show that $\ker \Pi(\alpha) = \bigwedge_k \ker \Pi(\alpha_k)$, that is, $= \bigcap_k \ker \Pi(\alpha_k)$. We have

$$\begin{aligned}
\xi \in \ker \Pi(\alpha) &\Leftrightarrow \xi \in \ker \Pi(\alpha)^{1/2} \\
&\Leftrightarrow \langle \Pi(\alpha)\xi, \xi \rangle = 0 \Leftrightarrow \sum_k \langle \Pi(\alpha_k)\xi, \xi \rangle = 0 \\
&\Leftrightarrow \langle \Pi(\alpha_k)\xi, \xi \rangle = 0, \ k \in \mathbb{N} \\
&\Leftrightarrow \xi \in \bigcap_k \ker \Pi(\alpha_k)^{1/2} = \bigcap_k \ker \Pi(\alpha_k).
\end{aligned}$$

Our use of $\mathbb{N}$ as an index set is no loss of generality since $\Pi(\alpha)$ depends on $\alpha$ only up to null sets and $\mu$ is $\sigma$-finite. Thus we have checked that $Y$ is the free projection-valued upper extremal process.

In view of the formula for the free Poisson distribution (Section 2.7 in [17]) we have that if $\mu(\alpha) = a < 1$, the spectrum of $\Pi(\alpha)$ is $\{0\} \cup [(1 - a^{1/2})^2, (1 + a^{1/2})^2]$ and $Y(\alpha)$ is the spectral projection of $\Pi(\alpha)$ for $[(1 - a^{1/2})^2, (1 + a^{1/2})^2]$. This gives $\|Y(\alpha) - \Pi(\alpha)\| \leq (1 + a^{1/2})^2 - 1 \leq 3a^{1/2}$ and $|Y(\alpha) - \Pi(\alpha)|_1 \leq \|Y(\alpha) - \Pi(\alpha)\|\tau(Y(\alpha)) \leq 3a^{1/2} \cdot a = 3a^{3/2}$.

To prove that the von Neumann algebras generated by $\{Y(\alpha) | \alpha \in \mathcal{B}\}$ and $\{\Pi(\alpha) | \alpha \in \mathcal{B}\}$ coincide, it suffices to show their $L^1$-spaces coincide, and hence it suffices to show $\Pi(\alpha)$ is in the $L^1$-closure of the linear span of the $Y(\alpha)$. Given $\alpha \in \mathcal{B}_f$, if $\mu$ is diffuse, there are pairwise disjoint $\alpha_1, \ldots, \alpha_k \in \mathcal{B}_f$ so that $\alpha_1 \cup \cdots \cup \alpha_k = \alpha$ and $\mu(\alpha_j) < \varepsilon$. Then we have

$$|\Pi(\alpha) - (Y(\alpha_1) + \cdots + Y(\alpha_k))|_1 \leq \sum_j |\Pi(\alpha_j) - Y(\alpha_j)|_1$$
$$\leq 3\varepsilon^{1/2}(\mu(\alpha_1) + \cdots + \mu(\alpha_k))$$
$$\leq 3\varepsilon^{1/2}\mu(\alpha),$$

which proves our assertion. $\square$

REMARK 5.3. Note that in the realization of the free Poisson process via circular elements, $Y(\alpha)$ being the range projection of $\Sigma c_\iota \chi(\alpha \cap \Omega_\iota) c_\iota^*$ is the same as being the range projection of $C\chi(\alpha)$, where $C$ is a column matrix with entries $C_\iota$ and $\chi(\alpha)$ is a diagonal matrix with entries $\chi(\alpha \cap \Omega_\iota)$ on the diagonal. This can also be translated into the asymptotic random matrix realization (see Section 7.3 in [17]): the $Y(\alpha)$ are the large $N$ limit of $\Gamma_N D_N(\alpha)$ [or equivalently of $\Gamma_N D_N(\alpha) \Gamma_N^*$], where $\Gamma_N$ is an appropriate $N \times [tN]$ matrix with i.i.d. Gaussian entries and $D_N(\alpha)$ are deterministic $[Nt] \times [Nt]$ projection matrices with joint limit distribution the same as the $\chi(\alpha)$'s in $L^\infty(\Omega, \Sigma, t^{-1}\nu)$ [here $t = \nu(\Omega) < \infty$; we leave to the reader the easy adaptation to $\nu(\Omega) = \infty$].

REMARK 5.4. By relaxing the condition $Y(\omega) \in \text{Proj}(M)$ to $Y(\omega) \in M_h$, we may consider more general free upper extremal process over a set, that is, in addition to the condition that $\omega \to Y(\omega)$ be an upper extremal process over $B/\sim$, we also require that, for pairwise disjoint $\{\omega_k\}_{k \in K}$, the $(Y(\omega_k))_{k \in K}$ be freely independent. Note that in view of Remark 4.4, $E(Y(\omega); (t, \infty))$ will then be a projection-valued upper extremal process indexed by $B/\sim$ and the freeness requirement for disjoint $\{\omega_k\}_{k \in K}$ also carries over to the $E(Y(\omega_k); (t, \infty))$.

An example of such a more general free upper extremal process over $(\mathcal{X}, \mathcal{B}, \mu)$ is the *free triangular upper extremal process* $Z(\omega)$ over $(\mathcal{X}, \mathcal{B}, \mu)$. It is characterized by the fact that $E(Z(\omega); (t, \infty)) = R_t(\omega)$ has the properties

$$t < 0 \Rightarrow R_t(\omega) = I,$$
$$t \geq 1 \Rightarrow R_t(\omega) = 0,$$
$$0 \leq t < 1 \Rightarrow \tau(R_t(\omega)) = \min((1-t)\mu(\omega), 1).$$

The kind of argument that we used for the existence of the projection-valued free process over $(\mathcal{X}, \mathcal{B}, \mu)$ also works for the triangular process (we leave the easy details to the reader).

Note also that, in addition to the triangular process, free upper extremal processes over $(\mathcal{X}, \mathcal{B}, \mu)$ can be obtained applying Remark 4.7 to $Z(\omega)$, that is, processes $f_*Z$, where $f$ is an increasing lower semicontinuous function $\mathbb{R} \to \mathbb{R}$.

**6. Free max-stable distributions and free max-domains of attraction.** It will be convenient to work with distribution functions and to adapt some of the notation and definitions of Section 3.

DEFINITION 6.1. If $F$ and $G$ are two distribution functions on the real line, we define the distribution function $F \boxvee G$ to be $(F + G - 1)_+$ and we define the $n$-fold iterate of this operation to be

$$F^{\boxvee n} = F \underbrace{\boxvee \cdots \boxvee}_{n \text{ times}} F = (nF - (n-1))_+.$$

We want to study the possible asymptotic behavior of $F^{\boxvee n}$ when $n$ tends to infinity. Let us start with some trivial properties of the tail and the support of $F^{\boxvee n}$. If $\bar{F} = 1 - F$ denotes the tail of the distribution function $F$, then

$$\overline{F \boxvee G} = (\bar{F} + \bar{G}) \wedge 1$$

and

$$\overline{F^{\boxvee n}} = n\bar{F} \wedge 1.$$

If $[\alpha(F), \omega(F)]$ denotes the support in $[-\infty, \infty]$ of the probability distribution defined by $F$, then for all $n \geq 2$ and $F$,

$$\alpha(F^{\boxvee n}) > -\infty$$

and

$$\omega(F^{\boxvee n}) = \omega(F).$$

More precisely,

$$\alpha(F^{\boxvee n}) = \sup\left\{x \in \mathbb{R} | F(x) \leq 1 - \frac{1}{n}\right\}$$

and

$$\lim_{n\to\infty} \alpha(F^{\boxtimes n}) = \omega(F).$$

In case there is $u_n$ such that

$$F(u_n) = 1 - n^{-1},$$

then

$$\alpha(F^{\boxtimes n}) \geq u_n$$

and one can easily interpret $F^{\boxtimes n}$ as a conditioned (or thresholded) distribution. Indeed, if $F$ is the distribution function of a random variable $X$, then

$$F^{\boxtimes n}(x) = P(X \leq x | X > u_n),$$

that is, $F^{\boxtimes n}$ is the distribution function of the random variable $X$ conditioned to be larger than the threshold $u_n$ (cf. Proposition 3.4).

DEFINITION 6.2. A distribution function $F$ is freely max-stable iff for every $n \geq 1$, there exist $a_n, b_n \in \mathbb{R}$, $a_n > 0$, such that

$$F^{\boxtimes n}(a_n x + b_n) = F(x).$$

REMARK 6.3. If $F$ is freely max-stable, then its support is bounded from below, since $\alpha(F^{\boxtimes 2}) > -\infty$.

DEFINITION 6.4. A distribution function $F$ is in the free max-domain of attraction of the distribution function $G$ if there exist $a_n, b_n \in \mathbb{R}$, $a_n > 0$, such that, as $n \to \infty$,

$$F^{\boxtimes n}(a_n x + b_n) \xrightarrow{w} G(x)$$

(i.e., convergence at every point of continuity of $G$). The free max-domain of attraction of $G$ will be denoted by $\text{Dom}_{\text{free}}(G)$.

THEOREM 6.5. *The following statements are equivalent:*

  (i) *$G$ is freely max-stable;*
 (ii) *$\text{Dom}_{\text{free}}(G) \neq \varnothing$;*
(iii) *$G \in \text{Dom}_{\text{free}}(G)$.*

PROOF. Clearly (i) $\Rightarrow$ (iii) $\Rightarrow$ (ii) and we are left with proving that (ii) $\Rightarrow$ (i). If $F \in \text{Dom}_{\text{free}}(G)$, then there are $a_n, b_n \in \mathbb{R}$, $a_n > 0$, so that at every continuity point $x$ of $G$, we have

$$\lim_{n\to\infty} n\bar{F}(a_n x + b_n) \wedge 1 = \bar{G}(x).$$

Hence, if $k \geq 1$ is an integer, then

$$\lim_{n \to \infty} n\bar{F}(a_n x + b_n) \wedge k^{-1} = \bar{G}(x) \wedge k^{-1}$$

and

$$\lim_{n \to \infty} n\bar{F}(a_{nk} x + b_{nk}) \wedge k^{-1} = k^{-1}\bar{G}(x).$$

Thus, by a slight generalization of Khintchine's law of types, which we give in Lemma 6.6, we infer the existence of $\alpha_k, \beta_k \in \mathbb{R}$, $\alpha_k > 0$, such that

$$\bar{G}(\alpha_k x + \beta_k) \wedge k^{-1} = k^{-1}\bar{G}(x)$$

or, equivalently,

$$k\bar{G}(\alpha_k x + \beta_k) \wedge 1 = \bar{G}(x),$$

which is the same as

$$G^{\boxtimes k}(\alpha_k x + \beta_k) = G(x).$$

So, modulo Lemma 6.6 we have proved $G$ is max-stable.

To formulate the extension of Khintchine's law of types, let us call $G$ a $c$-defective distribution function (here $0 < c < 1$) if $G$ is a nondecreasing, right continuous function on $\mathbb{R}$ so that

$$\lim_{x \to +\infty} G(x) = 1$$

and

$$\lim_{x \to -\infty} G(x) = c.$$

Then $G$ can be viewed a the distribution function of a random variable that takes the value $-\infty$ with probability $c$.

LEMMA 6.6. *If $F_n$ is a sequence of $c$-defective distribution functions and if $G$ and $G_*$ are $c$-defective, nondegenerate distribution functions such that, as $n \to \infty$, we have*

$$F_n(a_n x + b_n) \xrightarrow{w} G(x),$$
$$F_n(\alpha_n x + \beta_n) \xrightarrow{w} G_*(x)$$

*for real constants $a_n > 0$, $b_n, \alpha_n > 0$, $\beta_n$, then*

$$\lim_{n \to \infty} \frac{a_n}{\alpha_n} = a \quad and \quad \lim_{n \to \infty} \frac{\beta_n - b_n}{a_n} = b$$

*for some $a > 0$, $b \in \mathbb{R}$ and*

$$G_*(x) = G(ax + b).$$

This lemma is a trivial consequence of Khintchine's law of types. Indeed one reduces it to the usual result by defining

$$\tilde{G}(x) = \frac{G(x) - c}{1 - c} \quad \text{and} \quad \tilde{G}_*(x) = \frac{G_*(x) - c}{1 - c}.$$

To conclude the proof of the theorem, we apply the lemma to the sequence of $(1 - k^{-1})$-defective distribution functions $1 - (n\bar{F} \wedge k^{-1})$ and to the $(1 - k^{-1})$-defective distribution functions $1 - \bar{G} \wedge k^{-1}$ and $1 - \bar{G}k^{-1}$. □

DEFINITION 6.7. We will say that a distribution function $F$ is of free extreme-value type if $F$ is of the same type as one of the following classes of distributions:

*Type* I: The exponential distribution $F(x) = (1 - e^{-x})_+$.
*Type* II: The Pareto distribution $F(x) = (1 - x^{-\alpha})_+$ for some $\alpha > 0$.
*Type* III: The Beta law $F(x) = 1 - |x|^\alpha$ for $-1 \leq x \leq 0$ and some $\alpha > 0$.

THEOREM 6.8. *$G$ is freely max-stable iff $G$ is of free extreme-value type [i.e., there exist $a > 0$ and $b$ real constants so that $G(ax + b)$ is one of the distributions listed in Definition* 6.7*].*

PROOF. Each of the distributions in Definition 6.7 is freely max-stable by a straightforward computation. To prove the converse statement, we need the following lemma.

LEMMA 6.9. *If $G$ is freely max-stable, then there exist measurable functions $a(s) > 0$, $b(s)$, where $s \in [1, \infty)$, so that*

$$G^{\boxdot s}(a(s)x + b(s)) = G(x)$$

*with $G^{\boxdot s}$ defined for $s \in [1, \infty)$ by*

$$\overline{G^{\boxdot s}(y)} = s\overline{G(y)} \wedge 1.$$

PROOF. Whereas $G$ is freely max-stable, there exist real numbers $a_n > 0$ and $b_n$ such that

$$\bar{G}(x) = n\bar{G}(a_n x + b_n) \wedge 1.$$

If $u_s(x) = a_{[ns]}x + b_{[ns]}$, then

$$\bar{G}(x) = [ns]\bar{G}(u_s(x)) \wedge 1.$$

This easily gives

$$\lim_{n \to \infty} n\bar{G}(u_s(x)) \wedge s^{-1} = \bar{G}(x)s^{-1}.$$

Thus, we have
$$n\bar{G}(a_n x + b_n) \wedge s^{-1} = \bar{G}(x) \wedge s^{-1}$$
and
$$\lim_{n \to \infty} n\bar{G}(a_{[ns]}x + b_{[ns]}) \wedge s^{-1} = \bar{G}(x)s^{-1}.$$

By the extension of Khintchine's law of types in Lemma 6.6, we conclude that there exist measurable functions $a(s) > 0$ and $b(s)$ such that
$$s^{-1}\bar{G}(a(s)x + b(s)) = \bar{G}(x) \wedge s^{-1}$$
or, equivalently,
$$\bar{G}(a(s)x + b(s)) = s\bar{G}(x) \wedge 1,$$
that is,
$$G(a(s)x + b(s)) = G^{\boxplus s}(x). \qquad \square$$

By applying Lemma 6.9 to $G^{\boxplus st}$, one easily gets that
$$a(st) = a(t)a(s),$$
$$b(st) = a(t)b(s) + b(t)$$
if $s, t \in (1, \infty)$.

It is easily seen that if one extends $a(s)$ to $(0, \infty)$ by $a(1) = 1$ and $a(s) = (a(s^{-1}))^{-1}$ if $s < 1$, then the equation $a(st) = a(t)a(s)$ holds for all $s > 0$ and $t > 0$. Since $a(s)$ is measurable, it is known that there exists $\theta \in \mathbb{R}$ such that $a(s) = s^\theta$.

*Case* I. If $\theta = 0$, that is, $a(s) \equiv 1$, then $b(st) = b(s) + b(t)$ for $s, t > 1$ and $\tilde{b}(s) = e^{b(s)}$ satisfies $\tilde{b}(st) = \tilde{b}(s)\tilde{b}(t)$. Thus again, there exists $c \in \mathbb{R}$ such that $b(s) = -c \ln s$ and
$$G^{\boxplus s}(x) = G(x - \ln s)$$
for $s > 1$. It is easy to see that this implies that $G$ is of the same type as the exponential distribution. Indeed, we just check that $c > 0$, since
$$\bar{G}(x - c \ln s) = \overline{G^{\boxplus s}}(x) = s\bar{G}(x) \wedge 1 \geq \bar{G}(x),$$
which implies that $x - c \ln s \leq x$ and thus $c \geq 0$. If $G$ is nondegenerate, we cannot have $c = 0$. Then we see that $G(x) < 1$ for all $x$. Indeed, if $G(x) = 1$, then $\bar{G}(x) = 0$, so that $\overline{G^{\boxplus s}}(x) = 0$; that is, $G^{\boxplus s}(x) = 1$ for all $s \geq 1$. Hence $G(x - c \ln s) = 1$ for all $s \geq 1$, which gives $G(y) = 1$ for all $y \leq x$, which is not possible.

Now, if $a = \bar{G}(x) > 0$ and $y = x - c \ln s$, then
$$\bar{G}(y) = sa \wedge 1 = (ae^{x/c}e^{-y/c}) \wedge 1$$

so that $\alpha(G) = x + c \ln a$. This shows that

$$\ln \bar{G}(x) = (\alpha(G) - x)c^{-1} \quad \text{for } x \geq \alpha(G),$$

that is

$$\bar{G}(x) = e^{(\alpha(G)-x)c^{-1}}$$

and $G(\alpha(G) + cx)$ is exponentially distributed.

*Case* II. If $\theta > 0$, then $a(s) = s^\theta > 1$ for $s > 1$. From

$$b(st) = a(t)b(s) + b(t) = a(s)b(t) + b(s)$$

we see that $b(s)(1 - a(s))^{-1}$ is a constant $c$. We infer that $b(s) = c(1 - s^{-\theta})$ and $G^{\boxplus s}(x) = G(s^{-\theta}(x - c) + c)$. If $H(x) = G(x + c)$, then

$$H^{\boxplus s}(x) = H(s^{-\theta}x),$$

that is, $s\bar{H}(x) \wedge 1 = \bar{H}(s^{-\theta}x)$. If $x_0 < 0$ is such that $0 < \bar{H}(x_0) < 1$, then this shows that $\bar{H}(s^{-\theta}x_0)$ is an increasing function on some interval $(1, 1 + \varepsilon)$, since it is equal to $s\bar{H}(x_0)$ for $s$ close to 1. However, this function is nonincreasing if $x_0 < 0$. Thus for every $x_0 < 0$ we see that $\bar{H}(x_0) \in \{0, 1\}$. Since $H$ is nondegenerate, it is then impossible that $\bar{H}(0) = 0$. Choosing $s > 1$ such that $s^{-1} < \bar{H}(0)$, we have

$$1 = s\bar{H}(0) \wedge 1 = \bar{H}(s^{-\theta} \cdot 0) = \bar{H}(0).$$

Hence $\bar{H}(0) = 1$. This proves that $\alpha(H) \geq 0$.

Let $x_0 > 0$ be such that $0 < H(x_0) < 1$. Then if $y = s^{-\theta}x_0$, one sees that $\bar{H}(y) = s\bar{H}(x_0) \wedge 1$. Hence $\bar{H}(y) = 1$ iff $s \geq s_0 = (\bar{H}(x_0))^{-1}$, that is, iff $y \leq s_0^{-\theta}x_0$. This proves that $\alpha(H) > 0$ and $\alpha(H) = s_0^{-\theta}x_0$. This is valid for every $x_0$ such that $0 < H(x_0) < 1$, so that $\alpha(H) = (\overline{H(x)})^\theta x$ for $x \geq \alpha(H)$ and

$$\bar{H}(x) = (\alpha(H))^{1/\theta}x^{-1/\theta},$$

that is, $H$ is a Pareto distribution.

*Case* III. If $\theta < 0$, this is completely similar to Case II and $G$ is of the type of a Beta law. □

In particular we have also proved the following fact:

THEOREM 6.10 (Free extremal type theorem). *The following statements are equivalent:*

(i) *There exists a distribution $F$ and constants $a_n, b_n \in \mathbb{R}$, $a_n > 0$, such that $F^{\boxplus n}(a_n x + b_n) \xrightarrow{w} G(x)$ as $n \to \infty$.*

(ii) *Distribution $G$ is of the type of a free extreme-value distribution.*

This result in the free probability setting is the equivalent of the classical extremal type theorem. One can introduce a natural mapping that relates the two statements.

Given $c > 0$, we define a function on $[0, 1]$ by

$$f_c(u) = (1 + c \ln u)_+.$$

Then $f_c$ is nondecreasing and $f_c(u) = 0$ iff $u \leq e^{-1/c}$. Moreover $f_c(1) = 1$. We have

$$f_c(uv) = (1 + c \ln u + c \ln v)_+$$
$$= (f_c(u) + f_c(v) - 1)_+.$$

We can endow the set $[0, 1]$ with two semigroup structures: one that arises from usual multiplication and the other that arises from the operation related to free convolution:

$$u * v = (u + v - 1)_+.$$

Then $f_c$ (for any $c > 0$) is a homomorphism between these two semigroups on $(0, 1]$.

This homomorphism gives rise to a homomorphism between the semigroups of probability distribution functions endowed with either pointwise multiplication or the operation $\boxtimes$, which coincides with performing the operation $*$ pointwise. It follows that

$$f_c(FG) = f_c(F) \boxtimes f_c(G)$$

so that

$$f_c(F^n) = f_c(F)^{\boxtimes n}.$$

It is clear that the free extreme-value distribution functions are obtained from the classical ones by the map $f_1$. It is also clear that if $F$ is classically max-stable, then $f_c(F)$ is freely max-stable. If one could prove the converse directly, that is, that if $G$ is freely max-stable, it is the image by $f_c$ of a max-stable distribution-function $F$, then we could derive the free extremal type theorem directly from the classical one. We were able to obtain this result via a different route.

We will now prove that free max-domains of attraction and classical max-domains of attraction coincide for corresponding laws and that the normalizing constants are also equal.

We begin with Type I.

THEOREM 6.11. *The free max-domain of attraction of the exponential distribution coincides with the classical max-domain of attraction of the Gumbel distribution and the normalizing constants are equal.*

PROOF. $F$ is in the free max-domain of attraction of the exponential distribution iff there exist $b(s)$ and $a(s) > 0$ for $s > 1$ such that for all $x \geq 0$, we have

$$\lim_{s \to +\infty} s\bar{F}(b(s) + xa(s)) = e^{-x}.$$

(This is clear for $s$ running over $\mathbb{N}$ and we can take $a(s) = a([s])$, $b(s) = b([s])$.) In particular, $\lim_{s \to +\infty} s\bar{F}(b(s)) = 1$. Hence we may assume $b(s) < \omega(F)$ and that $b(s)$ is nondecreasing. It is easy to see that this implies that $F$ is in the domain of attraction of the exponential distribution iff there exists a function $g(t) > 0$ such that, for all $x > 0$, we have

$$\lim_{t \uparrow \omega(F)} \frac{\bar{F}(t + xg(t))}{\bar{F}(t)} = e^{-x}.$$

If $t < \omega(F)$, let us denote by $U$ the function $U(t) = 1/\bar{F}(t)$ so that

$$\lim_{t \uparrow \omega(F)} \frac{U(t + xg(t))}{U(t)} = e^x$$

and $U$ is $\Gamma$-varying (see [14], the definition in 0.4.3, page 26) on $(\alpha(F), \omega(F))$. By [14], Proposition 0.10, page 28, this is equivalent to $F$ being in the classical max-domain of attraction of the Gumbel law $\Lambda$. □

We now turn to the free max-domains of attraction of the Pareto distribution.

THEOREM 6.12. *The following statements are equivalent for $\alpha > 0$:*

(i) *$F$ is in the free max-domain of attraction of the Pareto distribution with exponent $\alpha$.*
(ii) *$F$ is in the max-domain of attraction of the Frechet distribution $\Phi_\alpha$.*
(iii) *$\bar{F}$ is $-\alpha$-regularly varying at $\infty$.*

*Moreover the normalization constants in* (i) *and* (ii) *can be chosen of the form $a_n = u_n$, $b_n = 0$.*

PROOF. The equivalence of (ii) and (iii) is a classical fact due to Gnedenko (see [14], the definition in 0.4.1 on page 13 and Proposition 1.11 on page 54).

To prove (iii) $\Rightarrow$ (i), assume $\bar{F}$ is $-\alpha$-regularly varying at $\infty$, that is,

$$\lim_{t \to \infty} \frac{\bar{F}(tx)}{\bar{F}(t)} = x^{-\alpha}$$

for all $x > 0$. Then $u_n = \inf\{t \in \mathbb{R} | \bar{F}(t) < n^{-1}\}$ will be such that

$$\lim_{n \to \infty} n\bar{F}(u_n) = 1$$

and hence
$$\lim_{n\to\infty} n\bar{F}(u_n x) = x^{-\alpha}$$
for all $x > 0$. It follows that
$$\lim_{n\to\infty} \overline{F^{\boxtimes n}}(u_n x) = \lim_{n\to\infty} \bar{F}(u_n x) \wedge 1 = x^{-\alpha} \wedge 1$$
if $x > 0$, and since $\overline{F^{\boxtimes n}}(u_n x)$ is decreasing and less than or equal to 1, the limit will be 1 if $x \in (-\infty, 1)$. This proves (i).

To conclude the proof we will show that (i) $\Rightarrow$ (iii). If for some choice of constants $a_n, b_n \in \mathbb{R}$, $a_n > 0$, we have
$$\lim_{n\to\infty} \overline{F^{\boxtimes n}}(a_n x + b_n) = x^{-\alpha} \wedge 1$$
for $x > 0$, then if $U = 1/\bar{F}$,
$$\lim_{n\to\infty} n^{-1} U(a_n x + b_n) = x^\alpha \qquad \text{when } x > 1.$$
Then if $V(y) = U^\leftarrow(y) = \inf\{s : U(s) \geq y\}$, we have
$$\lim_{n\to\infty} (V(y) - b_n) a_n^{-1} = y^{1/\alpha} \qquad \text{if } y > 1.$$
With $a(t) = a_{[t]}$ and $b(t) = b_{[t]}$, we then have
$$\lim_{t\to\infty} (V(ty) - b(t)) a(t)^{-1} = y^{1/\alpha}$$
and
$$\lim_{t\to\infty} (V(ty_1) - V(ty_2)) a(t)^{-1} = y_1^{1/\alpha} - y_2^{1/\alpha}$$
when $y, y_1, y_2 > 1$.

This implies
$$\lim_{t\to\infty} a(tx) a(t)^{-1}$$
$$= \lim_{t\to\infty} ((V(txy_1) - V(txy_2)) a(t)^{-1}) / (V(txy_1) - V(txy_2) a(tx)^{-1})$$
$$= ((xy_1)^{1/\alpha} - (xy_2)^{1/\alpha}) / (y_1^{1/\alpha} - y_2^{1/\alpha}) = x^{1/\alpha}.$$

Then, for any fixed $y_1 > y_2 > 1$, the function $V(ty_1) - V(ty_2)$ is $1/\alpha$-regularly varying at $\infty$ as a function of $t$. Since a function $W(t)$ is $1/\alpha$-regularly varying at $\infty$ iff $W(tz)$ for some fixed $z > 0$ is $1/\alpha$-regularly varying at $\infty$ as a function of $t$, we infer that $V(ty) - V(t)$ is $1/\alpha$-regularly varying at $\infty$ as a function of $t$ for all $y > 0$, $y \neq 1$. Thus we can use the last part of the proof of Proposition 1.11 in 1.2 of [14] (starting with the last two paragraphs of page 55 and continuing on pages 56 and 57) and conclude that $\bar{F}$ is regularly varying of index $-\alpha$ at $\infty$. $\square$

We finally turn to the domain of attraction of the Beta law (Type III).

THEOREM 6.13. *The following statements are equivalent:*

(i) $F$ *is in the free max-domain of attraction of* $G(x) = 1 - |x|^\alpha$ *for* $-1 \leq x \leq 0$ *and* $\alpha > 0$.

(ii) $F$ *is in the classical max-domain of attraction of the Weibull distribution* $\Psi_\alpha$.

(iii) $\omega(F) < \infty$ *and* $\bar{F}(\omega(F) - u)$ *is regularly varying of exponent* $\alpha$ *at* $0$, *that is,* $\lim_{h \downarrow 0} \frac{\bar{F}(\omega(F) - xh)}{\bar{F}(\omega(F) - h)} = x^\alpha$ *if* $x > 0$. *Moreover, the normalization constants can be chosen to be*

$$a_n = \omega(F) - u_n \quad \text{and} \quad b_n = \omega(F).$$

PROOF. The equivalence of (ii) and (iii) is a classical fact due to Gnedenko (see [14], Proposition 1.13 in 1.3).

To prove that (iii) implies (i), note that if $a_n = \omega(F) - u_n$ and $b_n = \omega(F)$, then

$$\overline{F^{\boxtimes n}}(a_n x + b_n) = n\bar{F}(\omega(F) + x(\omega(F) - u_n)) \wedge 1.$$

Since $u_n = \inf\{t : \bar{F}(t) < n^{-1}\}$, we get that if $h = \omega(F) - u_n$, then as $n \to \infty$ we have

$$\bar{F}(\omega(F) - h) = \bar{F}(u_n) \sim n^{-1}.$$

It follows that if $x < 0$, then

$$\lim_{n \to \infty} n\bar{F}(\omega(F) + x(\omega(F) - u_n)) = \lim_{h \to 0} (\bar{F}(\omega(F) - h))^{-1} \bar{F}(\omega(F) + xh) = |x|^\alpha.$$

Hence

$$\lim_{n \to \infty} \overline{F^{\boxtimes n}}(a_n x + b_n) = |x|^\alpha \wedge 1$$

if $x < 0$, which proves (i).

To prove the converse, that (i) implies (iii), we mimic the proof of the same statement in the Type II case, this time using the proof of Proposition 1.13 in [14], pages 59–62.

Assume that for some choice of $a_n, b_n \in \mathbb{R}$, $a_n > 0$, we have

$$\lim_{n \to \infty} F^{\boxtimes n}(a_n x + b_n) = G(x)$$

for $x \in \mathbb{R}$. Then if $U = 1/\bar{F}$ and $-1 \leq x \leq 0$, we have

$$\lim_{n \to \infty} n^{-1} U(a_n x + b_n) = |x|^{-\alpha}.$$

Hence if $V(y) = U^{\leftarrow}(y) = \inf\{s : U(s) \geq y\}$, we infer

$$\lim_{n \to \infty} (V(ny) - b_n) a_n^{-1} = -y^{-1/\alpha},$$

where $0 < y < 1$.

With $a(t) = a_{[t]}$ and $b(t) = b_{[t]}$, we have for $0 < y < 1$,

$$\lim_{t \to \infty} (V(ty) - b(t))(a(t))^{-1} = 1 - y^{-1/\alpha}$$

and for $0 < y_2 < y_1 < 1$,

$$\lim_{t \to \infty} (V(ty_1) - V(ty_2))(a(t))^{-1} = y_2^{-1/\alpha} - y_1^{-1/\alpha}.$$

This implies in turn that for $x > 0$, if $0 < xy_1 < 1$ and $0 < y_2 < y_1 < 1$, then

$$\lim_{t \to \infty} a(tx)(a(t))^{-1} = \lim_{t \to \infty} \frac{V(txy_1) - V(txy_2)}{a(t)} \lim_{t \to \infty} \frac{a(tx)}{V(txy_1) - V(txy_2)}$$
$$= ((xy_2)^{-1/\alpha} - (xy_1)^{-1/\alpha})(y_2^{-1/\alpha} - y_1^{-1/\alpha})^{-1} = x^{-1/\alpha}$$

and thus $a(t)$ is regularly varying at $\infty$ of index $-1/\alpha$.

Next we show that $\omega(F) < \infty$. Note that $\omega(F) = V(\infty) = \lim_{y \uparrow \infty} V(y)$. Since $a(\cdot)$ is regularly varying of index $-1/\alpha$, we infer that if $2^{-1/\alpha} < \lambda < 1$, then for some constant $A > 0$ we have

$$a(2^n) \leq A\lambda^n \qquad \text{for } n \in \mathbb{N}.$$

On the other hand, using

$$\lim_{n \to \infty} (V(2^{n+2}y_1) - V(2^{n+2}y_2))(a(2^n))^{-1} = y_2^{-1/\alpha} - y_1^{-1/\alpha}$$

with $y_1 = 1/2$ and $y_2 = 1/4$ gives that for some $B > 0$ and some $n_0 \in \mathbb{N}$, we have, for $n \geq n_0$,

$$V(2^{n+1}) - V(2^n) \leq Ba(2^{n+2}) \leq AB\lambda^{n+2} \qquad \text{for } n \in \mathbb{N}.$$

Clearly, since $0 < \lambda < 1$, this implies

$$\lim_{n \to \infty} V(2^n) - V(2^{n_0}) = \sum_{n \geq n_0} (V(2^{n+1}) - V(2^n))$$
$$\leq \sum_{n \geq n_0} AB\lambda^{n+2} < \infty,$$

that is, $\omega(F) < \infty$.

Since $V(ty_1) - V(ty_2)$, like $a(t)$, is regularly varying of index $-1/\alpha$, we infer as in the Type II case that $V(ty) - V(t)$ is also regularly varying of index $-1/\alpha$ for all $y > 0$, $y \neq 1$.

With these preparations the next step is to prove that $(\omega(F) - V(t))^{-1}$ is regularly varying of index $1/\alpha$ and thus that $\bar{F}(\omega(F) - s^{-1})$ as a function of $s$ is regularly varying of index $\alpha$ [this means $\bar{F}(\omega(F) - u)$ is regularly varying of exponent $\alpha$ at 0]. The proof can now be completed using the last part of the proof of Proposition 1.13 in 1.3 of [14], pages 61 and 62. $\square$

**7. Peaks over threshold.** The probability distributions we found as possible limits of free extremal convolutions are well known in statistics. They come under the name *generalized Pareto distributions* (see [7], Section 3.4, and [13], Section 1.4, for a textbook treatment of this subject). Statisticians have introduced a convenient parametrization of these distributions.

DEFINITION 7.1. For $\gamma \in \mathbb{R}$ the standard generalized Pareto distribution is defined by its distribution function $G_\gamma$ given by

$$G_\gamma(x) = 1 - (1 + \gamma x)^{-1/\gamma}$$

for

$$\begin{cases} x > 0, & \text{if } \gamma > 0, \\ 0 < x < |\gamma|^{-1}, & \text{if } \gamma < 0, \end{cases}$$

and if $\gamma = 0$, $G_0(x) = 1 - e^{-x}$ for $x > 0$.

These distributions appear as limits in the peaks-over-threshold approach to extreme-value theory, which we sketch briefly.

DEFINITION 7.2. If $X$ is a random variable with distribution function $F$, for $u < \omega(F)$ the exceedance (or excess) distribution function at threshold $u$ is given by

$$F^{[u]}(x) = P(X \leq u + x | X > u)$$
$$= \frac{F(u+x) - F(u)}{1 - F(u)}$$

for $x \geq 0$.

The main result of P.O.T. theory due to Balkema and De Haan [2]:
If $F^{[u]}(a_u x + b_u)$ *has a continuous limiting distribution function as* $u \uparrow \omega(F)$, *then*

$$\lim_{u \uparrow \omega(F)} \left| F^{[u]}(x) - G_\gamma\left(\frac{x}{\sigma_u}\right) \right| = 0$$

*for some shape and scale parameters $\gamma$ and $\sigma_u$.*

If $X_1, \ldots, X_n$ are $n$ i.i.d. random variables with common distribution function $F$ and if $X^{(1)} \geq X^{(2)} \geq \cdots \geq X^{(n)}$ are the order statistics, that is, the values of the variables ordered from largest to smallest, the traditional approach to extreme-value theory studies the distribution of the maximum $X^{(1)}$ or of the first $k$ maxima $X^{(1)}, \ldots, X^{(k)}$ properly normalized, when $n \to \infty$. The P.O.T. approach considers the distribution of the variables conditioned on being larger than a large threshold, which is very close to the free extremal convolution studied here (see the paragraph preceding Definition 6.2).

The only (technical) difference between the P.O.T. approach and the free probability extreme-value theory is that the latter does not fix a threshold $u$ for the random variable but rather for the value of the tail of the distribution function. A more serious difference is that the free approach introduces a binary operation on distribution functions and in some sense exhibits the P.O.T. theory as a result of the iteration of this operation.

The P.O.T. theory is very useful in various fields of statistics (insurance, reliability among others) and has been developed by R. L. Smith [15], A. C. Davison and R. L. Smith [6] and J. Pickands [12].

One should also notice that another (related) classical occurrence of the generalized Pareto distributions is as intensities of limiting Poisson point processes of extreme-value theory (see [14], page 210, Corollary 4.19). More precisely, $F$ is in the free domain of attraction of $G_\gamma$ iff the point measure

$$\sum_{i=1}^{n} \delta_{(X_i - a_n)/b_n}$$

converges to a Poisson point process with intensity $G_\gamma$.

**Acknowledgments.** D. Voiculescu did part of this work while holding an International Blaise Pascal Research Chair from the Ile de France region. Both authors were invited by A. Guionnet to a free probability event at ENS Lyon in June 2003, where this joint work was initiated.

## REFERENCES

[1] ANDO, T. (1989). Majorization, doubly stochastic matrices and comparison of eigenvalues. *Linear Algebra Appl.* **18** 163–248. MR0995373
[2] BALKEMA, A. A. and DE HAAN, L. (1974). Residual lifetime at great age. *Ann. Probab.* **2** 792–804. MR0359049
[3] BERCOVICI, H. and PATA, V. (with an appendix by P. Biane) (1999). Stable laws and domain of attraction in free probability theory. *Ann. of Math. (2)* **149** 1023–1060. MR1709310
[4] BERCOVICI, H. and VOICULESCU, D. (1992). Levy–Hinčin type theorems for multiplicative and additive free convolution. *Pacific J. Math.* **153** 217–248. MR1151559
[5] BERCOVICI, H. and VOICULESCU, D. (1993). Free convolution of measures with unbounded support. *Indiana Univ. Math. J.* **42** 733–773. MR1254116
[6] DAVISON, A. C. and SMITH, R. L. (1990). Models for exceedances over high thresholds (with discussion). *J. Roy. Statist. Soc. Ser. B* **52** 393–442. MR1086795
[7] EMBRECHTS, P., KLÜPPELBERG, C. and MIKOSCH, T. (1997). *Modelling Extremal Events for Insurance and Finance.* Springer, Berlin. MR1458613
[8] GUNAWARDENA, J., ED. (1998). *Idempotency.* Cambridge Univ. Press. MR1608365
[9] KADISON, R. V. and RINGROSE, J. (1986). *Fundamentals of the Theory of Operator Algebras.* Academic Press, Orlando, FL. MR0859186

[10] MASLOV, V. P. and SAMBORSKII, S. N., EDS. (1992). *Idempotent Analysis*. Amer. Math. Soc., Providence, RI. MR1203781

[11] NICA, A. and SPEICHER, R. (with an appendix by D. Voiculescu) (1996). On the multiplication of free $N$-tuples of noncommutative random variables. *Amer. J. Math.* **118** 799–837. MR1400060

[12] PICKANDS, J. (1975). Statistical inference using extreme value order statistics. *Ann. Statist.* **3** 119–131. MR0423667

[13] REISS, R. D. and THOMAS, M. (1997). *Statistical Analysis of Extremal Values*. Birkhäuser, Berlin. MR1464696

[14] RESNICK, S. I. (1987). *Extreme Values, Regular Variation and Point Processes*. Springer, New York. MR0900810

[15] SMITH, R. L. (1987). Estimating tails of probability distributions. *Ann. Statist.* **15** 1174–1207. MR0902252

[16] STRATILA, S. and ZSIDO, L. (1979). *Lectures on von Neumann Algebras*. Editura Academiei and Abacus Press. MR0526399

[17] VOICULESCU, D. (1998). *Lectures on Free Probability Theory. Lectures on Probability Theory and Statistics. Saint Flour XXVIII. Lecture Notes in Math.* **1783** 279–349. Springer, Berlin. MR1775641

[18] VOICULESCU, D., DYKEMA, K. and NICA, A. (1992). *Free Random Variables. CRM Monograph Series* **1**. Amer. Math. Soc., Providence, RI. MR1217253

COURANT INSTITUTE FOR MATHEMATICAL SCIENCES
NEW YORK UNIVERSITY
251 MERCER STREET
NEW YORK, NEW YORK 10012
USA
E-MAIL: benarous@cims.nyu.edu

DEPARTMENT OF MATHEMATICS
UNIVERSITY OF CALIFORNIA AT BERKELEY
BERKELEY, CALIFORNIA 94720-3840
USA
E-MAIL: dvv@math.berkeley.edu